\documentclass[11pt]{amsart}
\theoremstyle{plain}
\newtheorem{thm}{Theorem}[section]
\newtheorem{theorem}[thm]{Theorem}

\newtheorem{lemma}[thm]{Lemma}
\newtheorem{corollary}[thm]{Corollary}
\newtheorem{proposition}[thm]{Proposition}
\theoremstyle{definition}

\newtheorem{notation}[thm]{Notation}

\newtheorem{definition}[thm]{Definition}

\newtheorem{example}[thm]{Example}

\newtheorem{question}[thm]{Question}

\numberwithin{equation}{section}

\newcommand{\bi}{{\bf i}}
\newcommand{\bl}{{\bf l}}
\newcommand{\bn}{{\bf n}}
\newcommand{\bH}{{\bf H}}

\newcommand{\bP}{{\bf P}}
\newcommand{\bR}{{\bf R}}
\newcommand{\bX}{{\bf X}}
\newcommand{\bY}{{\bf Y}}

\newcommand{\sC}{{\mathcal C}}
\newcommand{\sD}{{\mathcal D}}

\newcommand{\sJ}{{\mathcal J}}
\newcommand{\sK}{{\mathcal K}}
\newcommand{\sL}{{\mathcal L}}

\newcommand{\sO}{{\mathcal O}}
\newcommand{\sP}{{\mathcal P}}

\newcommand{\sR}{{\mathcal R}}
\newcommand{\sS}{{\mathcal S}}
\newcommand{\sT}{{\mathcal T}}
\newcommand{\sU}{{\mathcal U}}
\newcommand{\sV}{{\mathcal V}}
\newcommand{\sW}{{\mathcal W}}
\newcommand{\sX}{{\mathcal X}}

\newcommand{\A}{{\mathbb A}}

\newcommand{\C}{{\mathbb C}}

\newcommand{\G}{{\mathbb G}}

\newcommand{\I}{{\mathbb I}}
\newcommand{\J}{{\mathbb J}}

\newcommand{\BP}{{\mathbb P}}

\newcommand{\T}{{\mathbb T}}
\newcommand{\U}{{\mathbb U}}
\newcommand{\V}{{\mathbb V}}
\newcommand{\W}{{\mathbb W}}
\newcommand{\X}{{\mathbb X}}

\newcommand{\fg}{{\mathfrak g}}

\def\Sym{\mathop{\rm Sym}\nolimits}

\def\Hom{\mathop{\rm Hom}\nolimits}

\title[Metabelian Lie groups and VMRT]{Partial compactification of metabelian Lie groups with prescribed varieties of minimal rational tangents}

\author{Jun-Muk Hwang}

\thanks{This work was supported by the Institute for Basic Science (IBS-R032-D1).}
\begin{document}

\begin{abstract} We study minimal rational curves on a complex manifold  that are tangent to a distribution. In this setting, the variety of minimal rational tangents (VMRT) has to be isotropic with respect to the Levi tensor of the distribution. Our main result is a converse of this: any smooth projective variety isotropic with respect to a vector-valued anti-symmetric form can be realized as VMRT of minimal rational curves tangent to a distribution on a complex manifold.  The complex manifold is constructed as a partial equivariant compactification of a metabelian group, which is a result of  independent interest.
\end{abstract}

\maketitle


\section{Introduction}

{\bf Convention} \begin{itemize} \item[1.] We work in the holomorphic category: all varieties are complex analytic and all maps are holomorphic, unless stated otherwise. Open sets refer to Euclidean topology and open sets in Zariski topology are called Zariski-open sets. A general point of a complex analytic set $X$ means a point in a  Zariski-open subset of $X$, i.e., the complement of a closed analytic subset in $X$.
\item[2.]
We use the following notation for projective space.  For a vector space $V$, its projectivization $\BP V$ is the set of 1-dimensional subspaces of $V$. If $\sV$ is a vector bundle on a complex manifold, the projectivization of each fiber of $\sV$ defines the projective bundle $\BP \sV.$ For a nonzero vector $w \in V$, the corresponding  point of $\BP V$ is denoted by $[w]$. For a subset $S \subset \BP V$, we write $S^+ \subset V \setminus \{ 0 \}$ for the tautological $\C^{\times}$-bundle $S^+ \to S$.   If $S$ is a submanifold in a neighborhood of a point $s \in S$, the affine tangent space $T_{w} S^+ \subset V$ of $S^+$ at a point $w \in S^+$  with $s = [w]$ is sometimes written as $T_{s} S^+ \subset V$.
\item[3.] For a  holomorphic map $f: X \to Y$ between complex manifolds,
its differential at $x\in X$ is denoted by ${\rm d}_x f : T_x X \to T_{f(x)} Y$. When ${\rm d}_x f$ is surjective for all $x \in X$, the vector subbundle ${\rm Ker}({\rm d} f) \subset T X$ is denoted by $T^f$.  \end{itemize}

\bigskip
Minimal rational curves and varieties of minimal rational tangents  have been  studied  usually on smooth projective varieties (see the survey papers  \cite{HM99}, \cite{Hw01} and \cite{Hw12}).
For the discussion in this paper, we need to generalize these concepts to
complex manifolds which are not necessarily projective algebraic varieties.
For this generalization, we  note first that the space ${\rm RatCurves}(X)$ of rational curves on a complex manifold $X$ is well-defined as an analytic space.  In fact, the construction of ${\rm RatCurves}(X)$ in Definition II.2.11 of \cite{Ko}   for an algebraic variety $X$ can be generalized to a complex manifold $X$  by using the  analytic versions of Hilbert scheme and Chow variety mentioned in Theorem I.5.7 and I.5.8 of \cite{Ko}. Using this, we introduce the following concepts.

\begin{definition}\label{d.mrc}
Let $X$ be a complex manifold of dimension $n >0$.    \begin{itemize}
\item[(i)] A rational curve $C \subset X$ is {\em unbendable} if its normalization map $f: \BP^1 \to X$ satisfies  $$f^* TX \cong \sO(2) \oplus \sO(1)^{\oplus d} \oplus \sO^{\oplus(n-d -1)}$$ for some nonnegative integer $d \leq n-1$.
\item[(ii)]
    An irreducible component $\sK$ of ${\rm RatCurves}(X)$ is a {\em family of unbendable rational curves} on $X$, if a general member of $\sK$ is an unbendable rational curve.
   \item[(iii)] In (ii), let $\sK_x \subset \sK, x \in X,$ be the complex analytic subset parameterizing
    members of $\sK$ passing through a point $x \in X$ and let $\sK_x^{\rm norm}$ be its normalization.   Let $\kappa \in \sK^{\rm norm}_x$ be a point corresponding to a  rational curve in $X$ through $x$ which is immersed at $x$. The immersed tangent direction at $x$ of this rational curve  is denoted by $\tau_x(\kappa) \in \BP T_x X.$  For a point $x\in X$, let  $\sC_x \subset \BP T_x X$ be the union of $\tau_x(\kappa)
    \in \BP T_x X$ for all such $ \kappa \in \sK^{\rm norm}_x .$  We call $\sC_x$  the {\em variety of minimal rational tangents} (abbr. VMRT) of $\sK$ at $x$. It is a locally closed subset in $\BP T_x X$.
      \item[(iv)] In (iii), let  $D$ be a  distribution on a Zariski-open subset $X^o \subset X,$ i.e., a vector subbundle $D \subset TX^o$.  If the VMRT  $\sC_x $ is contained in $\BP D_x$ for a general point $x \in X^o,$
 we say that $\sK$ is {\em  subordinate to the distribution $D$}.
\item[(v)]  An irreducible component $\sK$ of ${\rm RatCurves}(X)$ is a
 {\em family of minimal rational curves} on $X$, if in terms of  the associated universal family morphisms  $$ \sK \stackrel{ \rho_{\sK}}{\longleftarrow}  {\rm Univ}_{\sK} \stackrel{\mu_{\sK}}{\longrightarrow} X,$$
there exists an open subset $O_{\sK} \subset X$  such that $\mu_{\sK}^{-1}(O_{\sK})$ is nonempty and the restriction $\mu_{\sK}|_{\mu_{\sK}^{-1}(O_{\sK})} \to O_{\sK}$ is a proper holomorphic map.
    It is known that  a family of minimal rational curves is a family of unbendable rational curves (e.g. Corollary IV.2.9 of \cite{Ko} or Theorem 1.2 of \cite{Hw01}).
   \end{itemize}
\end{definition}

We are interested in families of minimal rational curves subordinate to a nontrivial distribution.
Our main interest is when the complex manifold $X$ is a projective algebraic variety, but
we need to consider (not necessarily algebraic) complex manifold  $X$ to use more flexible construction. When $\sK$ is a family of  minimal rational curves on a smooth projective variety $X$, the VMRT  $\sC_x \subset \BP T_x X$ at a general point $x \in X$ is a projective subvariety and its projective geometry  plays an important role in understanding the geometry of $X$ itself.
It is known (Corollary 1 of \cite{HM04}) that, when $X$ is a projective algebraic variety and $\sK$ is a family of minimal rational curves, the normalization of the VMRT at a general point of $X$ must be smooth and quite often, the VMRT itself is smooth. For example, this is the case when $X$ is a smooth projective subvariety of projective space $\BP^N$ and $\sK$ is a family of lines of $\BP^N$ covering $X$ (Proposition 1.5 of \cite{Hw01}). Thus it is reasonable to consider  the cases when the VMRT at a general point is smooth. Then a basic question is:

\begin{question}\label{q.basic} What kind of  smooth projective varieties can arise as the VMRT $\sC_x \subset \BP D_x$ in Definition \ref{d.mrc}, when $\sK$ is  of a family of minimal rational curves subordinate to a distribution $D $ on a Zariski-open subset in $X$?
\end{question}

First of all, what if the distribution is trivial, meaning $D = TX$? Then it is known that any smooth projective variety can be realized as VMRT. More precisely, for any smooth projective variety $S \subset \BP^{n-1}$, there exists a smooth projective variety $X^S$ of dimension $n$ with a family $\sK$ of minimal rational curves whose VMRT $\sC_x \subset \BP T_x X^S$ at a general point $x \in X^S$ is projectively isomorphic to $S \subset \BP^{n-1}$ (see Example \ref{e.flat} below). Thus the interesting case of   Question \ref{q.basic} is when  $D \neq TX$. Moreover, it  is reasonable to assume that $D$ is not integrable, i.e., does not come from a foliation on a Zariski-open subset $X^o \subset X$. For otherwise, we can reduce the study of  rational curves on $X$ belonging to $\sK$  to the study of rational curves lying on  lower-dimensional submanifolds of $X$,  i.e., the leaves of the foliation (see Theorem 6.7 of \cite{Hw12}).
When $D$ is not integrable,  the following  tensorial invariant associated to $D$ is nonzero.

\begin{definition}\label{d.Levi}
Let $D \subset TM$ be a distribution on a complex manifold $M.$
The Lie brackets of local sections of $D$ determine a homomorphism
$${\rm Levi}_x^D: \wedge^2 D_x \to T_x M/D_x,$$ called the {\em Levi tensor}  of $D$ at $x$.
By Frobenius' theorem,   the Levi tensor is not identically zero if and only if $D$ is not integrable. \end{definition}

A necessary condition for Question \ref{q.basic} in terms of the Levi tensor of the distribution is given by the following result.

\begin{theorem}\label{t.HM}
Let $X$ be a complex manifold and $D \subset TX^o$ be a  distribution   defined on a Zariski-open subset $X^o \subset X$. Let $\sK$ be a family of unbendable rational curves on $X$ subordinate to $D$.  Then for  a general point $x \in X$ and  any smooth point $\alpha \in \sC_x$, the affine tangent space $T_{\alpha} \sC_x^+ \subset D_x$ satisfies \begin{equation}\label{e.HM} {\rm Levi}_x^D(u,v) = 0 \mbox{ for any } u, v \in T_{\alpha} \sC_x^+ \subset D_x. \end{equation}
\end{theorem}

This is equivalent to  Proposition 2.4 of \cite{Hw01}, which is a reformulation of Proposition 9 of \cite{HM98}.
The discussion in \cite{Hw01} is under the assumption that $X$ is a projective algebraic variety and $\sK$ is a family of minimal rational curves, but the proof of Proposition 2.4 of \cite{Hw01}  works verbatim in the setting of Theorem \ref{t.HM}.
Note that Proposition 2.4 of \cite{Hw01} states
\begin{equation}\label{e.Hw}  {\rm Levi}_x^D(a,v) = 0 \mbox{ for any } a \in \alpha^{+} \mbox{ and }  v \in T_{\alpha} \sC_x^+ \subset D_x,\end{equation} instead of (\ref{e.HM}), but the proof there can be easily modified to give (\ref{e.HM}). Actually the two statements are equivalent.  To see that (\ref{e.Hw}) for any smooth point $\alpha \in \sC_x$ implies (\ref{e.HM})  for any smooth point $\alpha \in \sC_x$,  let $\Lambda \subset \wedge^2 D_x$ be the kernel of ${\rm Levi}_x^D$. Then (\ref{e.Hw}) implies
$a \wedge w \in \Lambda$ for any $ a \in \alpha^{+} $ and $  w $ in the image of the second fundamental form of $\sC_x$ at $\alpha$ (see the proof of Proposition 1.3.1 of \cite{HM99}). But given any $u, v  \in T_{\alpha} \sC_x^+$, we can find  arcs $a(t) \in \sC_x^+$ and $b(t) \in T_{a(t)} \sC_x^+$ parameterized by $t$ in the unit disc $ \Delta \subset \C$ such that $$ a(0) = a, \ \frac{{\rm d}}{{\rm d} t}|_{t=0} a(t) = u, \ b(0) =v  \mbox{ and }  \frac{{\rm d}}{{\rm d} t}|_{t=0} b(t) =w,$$ for some $w$ in the image of the second fundamental form. Since $a(t) \wedge b(t) \in \Lambda$ by (\ref{e.Hw}), we have $$\frac{{\rm d}}{{\rm d} t}|_{t=0} (a(t) \wedge b(t)) = u \wedge v + a \wedge w \  \in \Lambda.$$ This implies that $u \wedge v \in \Lambda$, proving (\ref{e.HM}).

  By Theorem \ref{t.HM},  (\ref{e.HM}) is a necessary condition for the submanifold $\sC_x$  in Question \ref{q.basic}.
The following is our main result, which says that (\ref{e.HM}) is actually a sufficient condition for Question \ref{q.basic}.

\begin{theorem}\label{t.analytic}
Let $\U$ and $\W$ be vector spaces and $\omega: \wedge^2 \W \to \U$ be a  homomorphism. Let  $S \subset \BP \W$ be a smooth projective variety satisfying at any point $\alpha \in S$, $$ \omega(u, v) = 0 \mbox{ for any } u, v \in T_{\alpha}S^+.$$
Then we can find a  complex manifold $\sX$ with a distribution $D \subset T\sX^o$ on a Zariski-open subset $\sX^0 \subset \sX$ and a family $\sK$ of minimal rational curves on $\sX$ subordinate to $D$ such that there exists an isomorphism $D_x \cong \W$ for each $x \in \sX^o$ under which the Levi tensor ${\rm Levi}^{D}_x: \wedge^2 D_x \to T_x \sX/D_x$ is isomorphic to $\omega$  and the variety of minimal rational tangents $\sC_x \subset \BP D_x$  is isomorphic to $S \subset \BP \W$. \end{theorem}

The complex manifold $\sX$ in Theorem \ref{t.analytic} is obtained as a partial  compactification of a metabelian Lie group (in the sense of Definition \ref{d.metabelian}) by compactifying certain affine curves on the Lie group which are tangent to a natural distribution on the Lie group. Such an equivariant partial compactification of a metabelian Lie group, explained in Theorem \ref{t.quotient}, has independent interest.
The followings are two  classical examples of this construction.

\begin{example}\label{e.flat}
This is Example 5.6 in \cite{Hw12}.
When $\U =0$, Theorem \ref{t.analytic} (and Theorem \ref{t.quotient}) can be applied to any smooth projective variety $S \subset \BP \W$.
Let $\BP^n, n= \dim \W, $ be the projectivization of $\W \oplus \C$ and regard $\BP \W$ as a hyperplane in $\BP^n$. Let $\pi: X^S \to \BP^n$ be the blowup of $S \subset \BP \W \subset \BP^n$ and let $H \subset X^S$ be the proper transform of the hyperplane $\BP \W $.
Let $\sK$ be the component of ${\rm RatCurves}(X^S)$ whose general members are the proper transforms of lines in $\BP^n$ intersecting $S$. Then $\sK$ is a family of minimal rational curves on $X^S$
whose VMRT at a general point is isomorphic to $S \subset \BP \W$. Examining the proof of Theorem \ref{t.analytic}, one can see that
the open subsets $\pi^{-1}(\BP^n \setminus \BP \W)  \subset X^S \setminus H$ correspond to $\sX^o \subset \sX$ in  Theorem \ref{t.analytic} applied to the case of $\U =0$.  \end{example}

\begin{example}\label{e.GP}
Let $G$ be a simple Lie group and let $P$ be a maximal parabolic subgroup associated to a long simple root of $G$ such that the radical $N$ of $P$ has depth 2, in other words, its Lie algebra $\bn$ satisfies $[\bn, \bn] \neq 0$ and $[\bn, [\bn, \bn]] =0$ (for example, an orthogonal Grassmannian $G/P$, other than hyperquadrics or spinor varieties, satisfies this condition). Then $N$ is a metabelian group in the sense of Definition \ref{d.metabelian}.  There exists a unique $G$-invariant distribution $D \subset T (G/P)$ with $ D \neq T (G/P)$.  Let $\sK$ be the space of lines on $G/P$ under the minimal projective embedding. Then $\sK$ is a family of minimal rational curves subordinate to $D$. The VMRT of $\sK$ gives a $G$-invariant fiber subbundle $\sC \subset \BP D$ whose fiber corresponds to the highest weight orbit of $P$ (Proposition 1 in \cite{HM02}). The action of $N$ on $G/P$ has an open orbit $\sX^o.$ Examining the proof of Theorem \ref{t.analytic}, one can see that the complex manifold $\sX$ in Theorem \ref{t.analytic} is the Zariski-open subset of $G/P$ which is the union of all lines intersecting $\sX^o$.
\end{example}

There are many other examples where Theorem \ref{t.analytic} can be applied. The following example shows that given any submanifold $Z \subset \BP V$, its third Veronese embedding $S = v_3(Z) \subset \BP \Sym^3 V$ is isotropic with respect to some nonzero vector-valued antisymmetric form on $\W = \Sym^3 V$ and Theorem \ref{t.analytic} can be applied.

\begin{example}\label{e.Veronese}
For a vector space $V$ of dimension $r \geq 2$, set $\W := \Sym^k V$ and let $S \subset \BP \W$ be the Veronese embedding of $\BP V$, the highest weight variety of the ${\rm SL}(V)$-representation on $\W$. Let $\W' \subset \wedge^2 \W$ be the subspace spanned by the planes tangent to $S^+ \subset \W$. Since ${\rm SL}(V)$ acts transitively on $\BP T (\BP V)$, the subspace $\W'$ must be an irreducible ${\rm SL}(V)$-module. Note that if $k \geq 3$, then $\wedge^2 \W = \wedge^2 (\Sym^k V)$ is not an irreducible ${\rm SL}(V)$-module, namely, we can write $$\wedge^2 \W = \W' \oplus \U$$ for some nonzero ${\rm SL}(V)$-module $\U$. Thus we have a nonzero projection $\omega: \wedge^2 \W \to \U$ such that $S$ is $\omega$-isotropic.
When $r=2$ and $k=3$, the $\omega$-isotropic submanifold $S$ is the twisted cubic curve and the resulting metabelian group  is equal to $N$ in Example \ref{e.GP} with $G$ a simple Lie group of type ${\rm G}_2$. If $r \geq 3$ or $k \geq 4$, the complex manifold $\sX$ in Theorem \ref{t.analytic} does not seem to have been studied. \end{example}

A special case of Theorem \ref{t.analytic} and Theorem \ref{t.quotient} when $D$ is a contact structure and $S$ is a Legendrian submanifold in $\BP \W$ is proved in \cite{Hw20} (see also \cite{HwMa} when $S$ is a Legendrian submanifold of special type). In Example \ref{e.GP}, this is the case when $G/P$ is a homogeneous contact manifold. The method in \cite{Hw20}, however,  works only for the contact distribution, because it employs the symplectic geometry of the cotangent bundle. For the general situation of Theorem \ref{t.analytic}, symplectic geometry cannot be used and Theorem \ref{t.quotient} requires a completely different argument which uses the tangent bundle instead of the cotangent bundle.

Theorem \ref{t.analytic} would give an affirmative answer to Problem 6.3 in \cite{Hw12} (also Question 1.2 in \cite{Hw20}), if the complex manifold $\sX$ is quasi-projective. Unfortunately, we do not know whether the complex manifold $\sX$ in Theorem \ref{t.analytic}  can be made quasi-projective.
All we know is that $\sX$ is very close to an algebraic variety in the sense that we can modify it to obtain a quasi-projective variety. This modification is used in the proof of the next theorem.

\begin{theorem}\label{t.algebraic}
Let $\U$ and $\W$ be vector spaces and $\omega: \wedge^2 \W \to \U$ be a surjective homomorphism. Let  $S \subset \BP W$ be a smooth projective variety satisfying at any point $\alpha \in S$, $$ \omega(u, v) = 0 \mbox{ for any } u, v \in T_{\alpha}S^+.$$
Then we can find a smooth projective variety $\bX$  with a distribution $D \subset T \bX^o$ on a Zariski-open subset $\bX^o \subset \bX$ and a family $\sK$ of unbendable rational curves subordinate to $D$ such that  there exists an isomorphism $D_x \cong \W$ for each $x \in \bX^o$ under which the Levi tensor ${\rm Levi}^{D}_x: \wedge^2 D_x \to T_x \bX/D_x$ is isomorphic to $\omega$  and the (closure of) variety of minimal rational tangents $\sC_x \subset \BP D_x$  is isomorphic to $S \subset \BP \W$.  \end{theorem}

Note that  the complex manifold $\bX$ is a projective algebraic variety, but $\sK$ is only a family of unbendable rational curves in Theorem \ref{t.algebraic}, while $\sK$ is a family of minimal rational curves in Theorem \ref{t.analytic}, but $\sX$ is not necessarily quasi-projective.

The content of the paper is as follows.  Basic results on horizontal lines in a metabelian Lie group are collected in Section \ref{s.MA}. Section \ref{s.analytic} studies the geometry of horizontal lines in the direction of an isotropic submanifold in the natural distribution on the metabelian Lie group.
In Section \ref{s.algebraic}, a partial compactification of a metabelian Lie group with a prescribed VMRT is explained in Theorem \ref{t.quotient} from which
the proofs of Theorem \ref{t.analytic} and Theorem \ref{t.algebraic} follow.

\section{Horizontal lines in metabelian Lie groups}\label{s.MA}
\begin{definition}\label{d.metabelian}
Let $\U$ and $\W$ be vector spaces and $\omega: \wedge^2 \W \to \U$ be a homomorphism. \begin{itemize} \item[(i)] Define a Lie algebra structure on $\fg := \W \oplus \U$ by  $$ [\fg, \U] =0 \mbox{ and } [w, w'] := \omega(w, w') \mbox{ for all } w, w' \in \W.$$ Following \cite{Ga}, we call $\fg$ a {\em metabelian Lie algebra}. \item[(ii)] We fix a simply connected Lie group $\G$ with Lie algebra  $\fg,$ called a {\em metabelian Lie group}. It is easy to see that the exponential map $\exp: \fg \to \G$ is a biregular isomorphism of algebraic varieties. Let $\log: \G \to \fg= \W \oplus \U$ be the inverse of the exponential map and write for $x \in \G$,   $$ \log (x) = x^W \oplus  x^U \mbox{ with } x^W \in \W \mbox{ and }x^U \in \U. $$ Then the group multiplication satisfies
$$\log (x_1 \cdot x_2) = (x_1^W + x_2^W) \oplus ( x_1^U + x_2^U + \frac{1}{2} \omega(x_1^W, x_2^W)).$$ It is convenient to write an element  $x \in \G$ as $(x^W, x^U)$ such that $\log (w,u) = w \oplus u$  for $w  \in \W$ and $u \in \U$. Then the group multiplication can be written as
$$(w, u) \cdot (w', u') = (w + w', u + u' + \frac{1}{2} \omega(w, w'))$$ for any $w, w' \in \W$ and $u, u' \in \U$. Denote by $o = (0,0) \in \G$  the identity element of $\G$.
\item[(iii)]   The vector subspace $\W \subset \fg =  T_o \G$ defines a left-invariant distribution $\sW \subset  T \G$. \end{itemize}
\end{definition}

\begin{definition}\label{d.line}
For a nonzero vector $w \in \W$, define the curve $\ell^w_o \subset \G$ as the underlying variety of the 1-parameter subgroup $\exp(\C w) \subset \G$: $$ \ell^w_o := \{ (tw, 0) \in \G \mid  t \in \C\}.$$ For any point $x = (x^W, x^U) \in  \G$, define $$\ell^w_x := x \cdot \exp(\C w) = \{(x^W + tw, x^U + \frac{t}{2} \omega(x^W,w) \mid  t \in \C\}.$$ A $w$-{\em line} on $\G$ means $\ell^w_x$ for some  $x \in\G$.
An affine curve on $\G$ is called a {\em horizontal line} if it is a $w$-line for some nonzero $w \in \W.$
\end{definition}

\begin{lemma}\label{l.W}
For $x = (x^W, x^U) \in \G$, the biregular isomorphism $\log: \G \to \fg$ induces an isomorphism $${\rm d}_x \log: T_x \G \to T_{\log (x)} \fg = \fg,$$  and the left translation $L_x: \G \to \G$ by $x$ gives the Maurer-Cartan isomorphism $${\rm d}_o L_x:  \fg = T_o \G  \to T_x \G.$$   Then for any vector $v \in \W \subset \fg$, we have $${\rm d}_x \log ({\rm d}_o L_x (v)) = v \oplus \frac{1}{2} \omega(x^W,  v) \ \in \fg.$$ \end{lemma}

\begin{proof} For $x = (x^W, x^U) \in \G,$ the map $\log \circ L_x: \G \to \fg$  sends the
 1-parameter group $\exp(t v) = \{(t v, 0) \in \G \mid t \in \C\}$ to
 $$ \{ (x^W + t v) \oplus  (x^U + \frac{t}{2} \omega(x^W, v)) \in \fg \mid t \in \C \}.$$
 By taking derivative with respect to $t $ at $t=0$, we obtain the result. \end{proof}

\begin{definition}\label{d.sL}
For each $s= [w] \in \BP \W$, let $\A_s= \A_{[w]} \subset \G$ be the algebraic subgroup $\exp (\C w)$ and let $\A \subset \BP \W \times \G$ be the disjoint union of such subgroups parameterized by $\BP \W$, namely, $$ \A := \{(s, g) \in \BP \W \times \G \mid  g \in \A_{s}\}.$$
We have the projection $\eta: \A \to \BP \W$ whose fiber $\eta^{-1}(s) \subset \A$ at  $s \in \BP \W$ is isomorphic to $\A_{s} \subset \G$ via the natural projection $\A \to \G$ .  Let $\pi: \sL \to \BP \W$ be the fiber bundle over $\BP \W$ whose  fiber $\pi^{-1}(s)$ at $s \in \BP \W$ is the coset space $\G/\A_{s}$.  In other words, we have a sequence of morphisms
$$ \begin{array}{ccccc}  \A & \subset & \BP \W \times \G &  \to & \sL \\ \eta \downarrow
& & \downarrow & & \downarrow \pi \\ \BP \W & =  & \BP \W & = & \BP \W\end{array} $$ which gives   the exact sequence $$0 \to \A_{s} \longrightarrow \G \longrightarrow \G/\A_{s} = \pi^{-1}(s) \to 0 $$ over each point $s\in \BP \W.$
By the left $\G$-action, there is a natural trivialization of the projective bundle on $\G$ \begin{equation}\label{e.PW}  \BP \sW \ \stackrel{\cong}{\longrightarrow} \ \BP \W \times \G. \end{equation} The pair of morphisms \begin{equation}\label{e.hline} \sL \stackrel{\varrho}{\longleftarrow} \BP \sW \stackrel{\upsilon}{\longrightarrow} \G,\end{equation} where $ \upsilon$ is the natural projection in (\ref{e.PW}) and
 $\varrho$ is the composition of the projection $\BP \W \times \G \to \sL$ and (\ref{e.PW}),  can be viewed as the universal family of horizontal lines on $\G$. In other words, we can regard  $\sL$ as the set of all horizontal lines on $\G$ such that the fiber $\pi^{-1}([w])$ of $\pi: \sL \to \BP \W$ at $[w] \in \BP \W$ can be viewed as  the set of all $w$-lines  on $\G$.
\end{definition}

\begin{definition}\label{d.isotropic}
Let $S \subset \BP \W$ be a locally closed submanifold. \begin{itemize}
  \item[(i)] An $S$-{\em line} on $\G$ means a $w$-line for some $w \in S^+$.
\item[(ii)] The tangent spaces of $S$-lines form a fiber subbundle $\sS \subset \BP \sW.$
Restricting (\ref{e.hline}) to $\sS$, we obtain
$$ \begin{array}{ccccccc} S & \stackrel{\lambda}{\longleftarrow} &
\sR & \stackrel{\rho}{\longleftarrow} & \sS & \stackrel{ \mu}{\longrightarrow} & \G \\
\cap & & \cap & & \cap & & \| \\ \BP \W & \stackrel{\pi}{\longleftarrow} &  \sL & \stackrel{\varrho}{\longleftarrow} & \BP \sW & \stackrel{\upsilon}{\longrightarrow} & \G, \end{array} $$ where $\sR $ is   the submanifold $\varrho(\sS)$ of $\sL$.
The pair of morphisms  $$ \sR \ \stackrel{\rho}{\longleftarrow} \ \sS \ \stackrel{ \mu}{\longrightarrow} \G$$ can be viewed as the universal family of $S$-lines on $\G$.
\item[(iii)]  Denote by $$\iota: \sS\to S \times \G$$  the biregular isomorphism induced by (\ref{e.PW}).
 \end{itemize} \end{definition}

The following lemma is straightforward.

\begin{lemma}\label{l.identity} In Definition \ref{d.isotropic},
for  a point $\alpha \in \sS$  with  $\iota(\alpha) = (s, x) \in S \times \G,$ denote by $\ell \in \sR$ the point $\rho(\alpha)$ corresponding to the line $\ell^w_x$ with $s=[w] \in S$. \begin{itemize} \item[(i)] The endomorphism of $T_s S$ induced by the derivative ${\rm d}_{\alpha} \rho $, $$T_s S \subset (T_s S \oplus T_x \G)  \stackrel{{\rm d} \iota}{=}  T_{\alpha} \sS \stackrel{{\rm d}_{\alpha} \rho }{\longrightarrow}   T_{\rho(\alpha)} \sR  \stackrel{{\rm d} \lambda}{\longrightarrow} T_s S$$ is the identity of $T_s S.$
\item[(ii)] The inverse of   ${\rm d}_o L_x: \fg \to T_x \G$ induces an isomorphism
$$\gamma_{\ell}: {\rm Ker}({\rm d}_{\ell} \lambda: T_{\ell} \sR \to T_{\lambda(\ell)} S) = T_x \G/T_x(x \cdot \A_w)  \stackrel{({\rm d}_o L_x)^{-1}}{\longrightarrow}  \fg/\C w.$$
\item[(iii)] The homomorphism ${\rm d}_{\alpha} \rho$ sends the  subspace  $T_x \G$ of $T_{\alpha} \sS \stackrel{{\rm d}_{\alpha} \iota}{=} T_s S \oplus T_x \G$ onto ${\rm Ker}({\rm d}_{\ell} \lambda)$ and the composition
$$\fg \stackrel{{\rm d}_o L_x}{\longrightarrow} T_x \G \stackrel{{\rm d}_{\alpha} \rho \circ ({\rm d}_{\alpha}\iota)^{-1}}{\longrightarrow} T_x \G/T_x(x \cdot \A_w) \stackrel{\gamma_{\ell}}{\longrightarrow} \fg/\C w$$
is just the quotient homomorphism $\fg \to \fg/\C w.$ \end{itemize}
\end{lemma}

\section{Lines in the direction of an isotropic submanifold}\label{s.analytic}
\begin{definition}\label{d.S} In Definition \ref{d.metabelian}, a locally closed submanifold  $S \subset \BP \W$ is  {\em $\omega$-isotropic}   if $\omega(u,v) =0$ for any $s \in S$ and any $u, v \in T_{s} S^+ \subset \W$. \end{definition}

\begin{proposition}\label{p.h_t} In Definition \ref{d.isotropic}, fix an $\omega$-isotropic submanifold $S \subset \BP \W$. For any  $x \in \G$, $w \in S^+$ with $s= [w] \in S$, $v \in T_s S$ and $ t\in \C$,
define \begin{itemize} \item[(1)] $x_t := x \cdot \exp(t w) \in \G$ such that $\ell^w_x = \{ x_t \mid t \in \C\}$; \item[(2)] $\alpha_t := [T_{x_t} \ell^w_x] \in \sS;$  \item[(3)]
   $v_t \in T_{\alpha_t} \sS$ as the unique vector satisfying ${\rm d}_{\alpha_t} \iota (v_t) = v \oplus 0 \in T_s S \oplus T_{x_t} \G;$ and \item[(4)] $ h_t (v) := {\rm d}_{\alpha_t} \rho (v_t) - {\rm d}_{\alpha_0} \rho(v_0)  \in T_{\ell} \sR,  $ where $\ell = \rho(\alpha_t) = \rho(\alpha_0) \in \sR$  corresponds to the line $\ell^w_x$. \end{itemize} Here, $h_t$  determines  a homomorphism $ h_t: T_s S \to {\rm Ker}({\rm d}_{\ell} \lambda)$
by lemma \ref{l.identity} (i).  Then for any $v \in T_s S$ and $ t \in \C$, we have
 $$  \gamma_{\ell} \circ h_t (v) =  - t v (w) \ \in T_w S^+/\C w  \subset \fg/\C w,$$  where $\gamma_{\ell}$ is the isomorphism in Lemma \ref{l.identity} (ii) and $v(w) \in T_w S^+/\C w \subset  \fg/\C w$ denotes the evaluation of $$v \in T_s S = \Hom( \C w, T_w S^+/\C w) \subset \Hom(\C w, \fg/\C w)$$ at $w \in \C w$.
 \end{proposition}

  \begin{proof}
To prove the proposition, we may replace $S$ by a neighborhood of $s$ in $S$.
Choose a local section $\sigma: S \to S^+ \subset \W$ of  the $\C^{\times}$-bundle $S^+ \to S$ such that $\sigma(s) = w \in S^+.$ To simplify the notation, we identify $S$ with the submanifold  $\sigma(S) \subset S^+ \subset \W$ passing through $s =w$. Then $T_s S$ is identified with the subspace $T_w (\sigma(S)) \subset \W$ such that $v \in T_s S$ is identified with a vector in $\W$ which is sent to  $v(w) \in T_w S^+/\C w$ modulo $\C w.$ It suffices to prove \begin{equation}\label{e.prove}
\gamma_{\ell} \circ h_t (v) = -t v \mod  \C w. \end{equation}

For any $y \in S$ and $z \in \G$, the tangent vector to the line $\ell^y_z = \{ z \cdot \exp (ty) \mid t \in \C\}$ at $t = t_0$ is just the left translation of $y \in \fg$ by $z \cdot \exp (t_0 y) \in \G$. Thus the curve in $\sS \subset \BP T \G$ given by the tangent vectors to $\ell^y_z$ is sent by $\iota: \sS \cong S \times \G$ to
$$\{ (y, z \cdot \exp (t y)) \mid t \in \C\} \ \subset S \times \G.$$
Let us consider the action $\phi: \C \times \sS \to \sS$ of the additive group $\C$ on $\sS$ such that the action of $t \in \C$ on $ S \times \G$ induced by the isomorphism $\iota: \sS \cong S \times \G$  sends $(y, z) \in S \times \G$ with $ z = (z^W, z^U)$ to
\begin{equation}\label{e.phi} \phi_{t} (y,z) := (y, z \cdot \exp (ty)) = (y, (z^W + ty, z^U + \frac{t}{2} \omega(z^W, y))). \end{equation}
Then the $\phi$-orbit through $\iota^{-1} (y,z) \in \sS$   is the curve in $\sS$ given by the tangent vectors to the line $\ell^y_z$. Thus the morphism $\rho: \sS  \to \sR$ is the quotient by this $\C$-action $\phi.$
It follows that  $\alpha_t = \phi_t(\alpha_0)$ and \begin{equation}\label{e.quot}
{\rm d}_{\alpha_0} \rho (v_0) = {\rm d}_{\alpha_t} \rho ({\rm d}_{\alpha_0} \phi_t (v_0)) \end{equation} for any $ t \in \C$.

   Let $\Delta $ be the unit disc in $\C$ and choose an arc $$\{ w_{\tau} := w + \tau v +  (\cdots) \mid \tau \in \Delta\} \subset S \subset  \W,$$ where $(\cdots)$ stands for terms involving $\tau^k, k \geq 2$. This is an arc
   on $S \subset S^+$ passing through $w$ in the direction of $v \in T_s S$.
   Then $$\{(w_{\tau}, x) \mid \tau \in \Delta\} \subset  S \times \{x\} \subset S \times \G $$ can be regarded as an arc in $\mu^{-1}(x) \subset \sS$  by $\iota : \sS \cong S \times \G$. By (\ref{e.phi}), the action $\phi_t, t \in \C,$ sends this arc to the arc
$$  \{ (w_{\tau}, (x^W+ t w_{\tau}, x^U + \frac{t}{2} \omega(x^W, w_{\tau})) \mid \tau \in \Delta \} \subset S \times \G $$ passing through $(w, x_t)$.
Under the embedding into the vector space
$ S \times \G \stackrel{\sigma \times \log}{\longrightarrow} \W \oplus \fg,$ this arc is sent to
$$\{ w_{\tau} \oplus ((x^W+ t w_{\tau}) \oplus (x^U + \frac{t}{2} \omega(x^W, w_{\tau}))) \in \W \oplus \fg \mid \tau \in \Delta \}.$$
  Taking derivative with respect to $\tau$ at $\tau =0$, we obtain
$$v \oplus t (v \oplus \frac{1}{2} \omega (x^W, v)) \in \W \oplus \fg.  $$ This represents the image of
$ {\rm d}_{\alpha_0}\phi_t (v_0) \in T_{\alpha_t} \sS$ under the composition of two isomorphisms $$T_{\alpha_t} \sS \stackrel{{\rm d} \iota}{\longrightarrow} T_s S \oplus T_{x_t} \G \stackrel{{\rm d} ({\rm Id}_S \times \log)}{\longrightarrow} T_s S \oplus \fg,$$ in other words, \begin{equation}\label{e.arc} {\rm d}_{(w, x_t)} ({\rm Id}_S \times \log) \circ  {\rm d}_{\alpha_t}  \iota \circ {\rm d}_{\alpha_0} \phi_t (v_0) = v \oplus t (v \oplus \frac{1}{2} \omega(x^W, v)). \end{equation}

Set $v^t := {\rm d}_o L_{x_t}(v) \in T_{x_t} \G$ such that
\begin{equation}\label{e.v^t}  {\rm d}_{x_t} \log (v^t) \ = \ {\rm d}_{x_t} \log ({\rm d}_o L_{x_t} (v))  \ = \  v \oplus \frac{1}{2} \omega(x_t^W, v) \ \in \fg,\end{equation} where the second equality is from Lemma \ref{l.W}.
As $S$ is $\omega$-isotropic, we have $\omega(w, v) =0$ and
$$ \omega( x^W_t, v) = \omega(x^W + tw, v) = \omega(x^W,v).$$
Thus (\ref{e.arc}) and (\ref{e.v^t}) give \begin{eqnarray*} {\rm d}_{(w,x_t)} ({\rm Id}_S \times \log) \circ {\rm d}_{\alpha_t} \iota \circ  {\rm d}_{\alpha_0} \phi_t(v_0) & = & v \oplus t (v
\oplus \frac{1}{2} \omega(x_r^W, v)) \\ & = & v \oplus {\rm d}_{x_t} \log (t v^t) \\ &= &
{\rm d}_{(w, x_t)}({\rm Id}_S \times \log) ( v \oplus t v^t).\end{eqnarray*} Consequently, we have \begin{eqnarray*} {\rm d}_{\alpha_0} \phi_t(v_0) &=& ({\rm d}_{\alpha_t} \iota)^{-1} (v \oplus tv^t) \\ &=& ({\rm d}_{\alpha_t} \iota)^{-1}(v \oplus 0) + t ({\rm d}_{\alpha_t} \iota)^{-1}(0 \oplus v^t)  \\ &=&  v_t + t ({\rm d}_{\alpha_t} \iota)^{-1}(0 \oplus v^t) \ \in T_{\alpha_t} \sS,\end{eqnarray*} where we use  ${\rm d}_{\alpha_t} \iota (v_t) = v \oplus 0 \ \in \  T_s S \oplus T_{x_t} \G$ from our choice of $v_t$.
By (\ref{e.quot}),
\begin{eqnarray*} {\rm d}_{\alpha_0} \rho(v_0) &=& {\rm d}_{\alpha_t}\rho ({\rm d}_{\alpha_0} \phi_t (v_0)) \\ & = & {\rm d}_{\alpha_t} \rho (v_t + t ({\rm d}_{\alpha_t} \iota)^{-1}(0 \oplus v^t)) \\ & = & {\rm d}_{\alpha_t} \rho (v_t) + t {\rm d}_{\alpha_t} \rho (({\rm d}_{\alpha_t} \iota)^{-1}(0 \oplus v^t)), \end{eqnarray*}
from which we obtain $$h_t(v) = {\rm d}_{\alpha_t} \rho (v_t) - {\rm d}_{\alpha_0} \rho (v_0) = - t {\rm d}_{\alpha_t} \rho \circ  ({\rm d}_{\alpha_t} \iota)^{-1}(0 \oplus v^t). $$
It follows that  \begin{eqnarray*} \gamma_{\ell} (h_t (v)) &=& - t \gamma_{\ell} ({\rm d}_{\alpha_t} \rho \circ ({\rm d}_{\alpha_t} \iota)^{-1}(0 \oplus v^t))) \\ & = & - t \gamma_{\ell}({\rm d}_{\alpha_t} \rho \circ ({\rm d}_{\alpha_t} \iota)^{-1} (0 \oplus {\rm d}_o L_{x_t} (v))) \\& = &  - t v \mod \C w, \end{eqnarray*}
where the last equality is by Lemma \ref{l.identity} (iii).  This proves (\ref{e.prove}).
\end{proof}

\begin{definition}\label{d.sT}
Fix an $\omega$-isotropic submanifold $S \subset \BP \W$ of dimension $d >0$.
\begin{itemize} \item[(i)] For each $s \in S,$ the subspace $T_s S^+ \subset \W$ is an abelian Lie algebra of $\fg.$ Define $\T_s := \exp(T_s S^+),$   which  is an abelian algebraic subgroup of $\G$.
\item[(ii)] The cosets of $\T_s$  define a foliation of rank $d +1$ on $\G$ and as $s$ varies in $ S$, the cosets define a foliation of rank $d+1$ on $S \times \G$. Since $\A_s \subset \T_s$, this foliation descends under $\rho \circ \iota^{-1}$ of Definition \ref{d.isotropic} to  a distribution $$\sT \subset  {\rm Ker}({\rm d} \lambda) \subset T\sR$$ which defines a foliation of rank $d $ on $\sR.$ \end{itemize} \end{definition}

The following is immediate from Proposition \ref{p.h_t}.

    \begin{corollary}\label{c.sT}
    In Definition \ref{d.sT}, let $h_t$ be as in Proposition \ref{p.h_t}. Then  $h_t(T_s S) = \sT_{\ell} \subset {\rm Ker}(d_{\ell} \lambda)$ for any $s \in S$ and $\ell \in \lambda^{-1}(s) \subset \sR$. \end{corollary}

We skip the proof of the following easy lemma.

\begin{lemma}\label{l.tensor}
Let $V$ and $\widetilde{V}$ be vector spaces with a surjective homomorphism $\varphi: \widetilde{V} \to V.$ Suppose that we are given an injective homomorphism $j_{\infty}: V \to {\rm Ker}(\varphi)$ and a  holomorphic family of homomorphisms $\{ j_t: V \to \widetilde{V} \mid t \in \C\}$ such that   $$\varphi \circ j_t = {\rm Id}_V \mbox{ and }  j_t(v) - j_0(v) = t j_{\infty}(v)$$  for each $ t \in \C$ and $ v \in V$. Then we can find a  vector space $E= \C e_1 + \C e_2$ with a basis $\{ e_1, e_2\}$ such that the subspace $j_0(V) +  j_{\infty}(V)$ of $\widetilde{V}$ admits a tensor decomposition
$$j_0(V) +  j_{\infty}(V) = j_{\infty}(V) \otimes E$$ with $j_{\infty} (v) = v \otimes e_2$ and $ j_t(v) = v \otimes (e_1 + t e_2)$ for each $t \in \C$. \end{lemma}

\begin{corollary}\label{c.tensor} Let $E= \C e_1 + \C e_2$ be a vector space with a basis $\{e_1, e_2\}$.  In the setting of  Corollary \ref{c.sT},
there exists a distribution $\sD \subset T \sR$ of rank $2d,  d = \dim S,$  such that \begin{itemize} \item[(i)]
$\sD \cap {\rm Ker}(d \lambda) = \sT,$ and  \item[(ii)]
there is a tensor decomposition
$\sD_{\ell} \cong \sT_{\ell} \otimes E$ for each $\ell \in \sR$ which sends $\sT_{\ell} \subset \sD_{\ell}$  to $\sT_{\ell} \otimes e_2$ and   ${\rm d}_{\alpha_t} \rho (T_s S)$ to
$\sT_{\ell} \otimes (e_1 + t e_2)$ for each $t \in \C$. \end{itemize} \end{corollary}

\begin{proof}
We can apply Lemma \ref{l.tensor} to Proposition \ref{p.h_t} by setting
$$V = T_s S, \ \widetilde{V} = T_{\ell} \sR, \ \varphi = {\rm d}_{\ell} \lambda, \   j_t = {\rm d}_{\alpha_t} \rho|_{T_s S} $$ and
 $$j_{\infty}: T_s S = \Hom(\C w, T_s S^+/\C w) \stackrel{{\rm ev}_w}{\longrightarrow} T_s S^+/\C w \stackrel{\gamma_{\ell}^{-1}}{\longrightarrow} \sT_{\ell} \subset \widetilde{V}, $$ where ${\rm ev}_w$ is the evaluation at $w$ and $\gamma_{\ell}$ is from Lemma \ref{l.identity}. Then define $\sD_{\ell} $ as $j_0(V) + j_{\infty}(V).$
\end{proof}

\begin{definition}\label{d.sP} In the setting of Corollary \ref{c.tensor}, let ${\rm Gr}(d, T\sR)$ be the Grassmannian bundle of $d$-dimensional subspaces in $T \sR$.
\begin{itemize} \item[(i)] Corollary \ref{c.tensor} determines  a $\BP^1$-bundle $\sP \subset {\rm Gr}(d, T\sR)$ with the natural projection  $\psi: \sP \to \sR$ and  a distinguished section $\sigma: \sR \to \sP$  such that  for each $\ell \in \sR$, \begin{itemize}
\item[(1)] the point
$\sigma(\ell) $ corresponds to the $d$-dimensional subspace $\sT_{\ell} \subset T_{\ell} \sR$  and  \item[(2)] the complement of
$ \sigma(\ell)$ in $\sP_{\ell} := \psi^{-1}(\ell)$ corresponds to the images of $T_s S$ under ${\rm d} \rho$
along $\rho^{-1}(\ell).$  \end{itemize} Denote  by $\Sigma$ the image $\sigma (\sR) \subset \sP$. By abuse of notation, we denote by $\sT \subset T\Sigma$ the image ${\rm d} \sigma (\sT)$ of the distribution $\sT \subset T \sR$.   \item[(ii)] Define  a holomorphic map $\chi: \sS \to \sP \setminus \Sigma$
by $$\chi(\alpha) = [{\rm d}_{\alpha} \rho ( {\rm Ker}({\rm d}_{\alpha} \mu))] = [{\rm d}_{\alpha} \rho (T_s S)]  \in {\rm Gr}(d, T \sR) $$ for each $ \alpha  \in \sS$ with $s = \lambda (\rho(\alpha)).$ Corollary \ref{c.tensor} shows that $\chi$ is a biholomorphic map. \item[(iii)] Define a subbundle  $\sV^o \subset T(\sP \setminus \Sigma)$ of rank $d$ as ${\rm d} \chi (T^{\mu})$ where $T^{\mu} = {\rm ker}({\rm d} \mu) \subset T \sS$. Note that $\sV^o$ defines a foliation on $\sP \setminus \Sigma$ whose leaves are the $\chi$-images of the fibers of $\mu: \sS \to \G$. \end{itemize}
\end{definition}

\begin{lemma}\label{l.tautological}
In the setting of Lemma \ref{l.tensor}, let ${\rm Gr}(d, \widetilde{V}), d = \dim V,$ be the  Grassmannian of $d$-dimensional subspaces in $\widetilde{V}$
and let $\sU$ be the tautological vector bundle on ${\rm Gr}(d, \widetilde{V})$, which is  a subbundle of rank $d$ in the trivial bundle $\widetilde{V} \times {\rm Gr}(d, \widetilde{V}).$
Consider the holomorphic map $j: \BP E \to  {\rm Gr}(d, \widetilde{V})$ given by
$$j([e_2]) = [j_{\infty}(V)] \in {\rm Gr}(d, \widetilde{V}) \mbox{ and } j([e_1 + t e_2]) = [j_t(V)] \in {\rm Gr}(d, \widetilde{V})$$ for $t \in \C$. Then $j^* \sU \cong \sO(-1)^{\oplus d}$ on $\BP^1 = \BP E$. \end{lemma}

\begin{proof}
Fix a basis $\{ v_1, \ldots, v_d\}$ of $V$. For each $1 \leq i \leq d,$  the family of vectors $\{v_i \otimes (e_1 + t e_2) \mid t \in \C\} $ together with $v_i \otimes e_2$ span a line subbundle of $j^* \sU$ isomorphic to the tautological line bundle on $\BP E$. Thus $j^* \sU \cong \sO(-1)^{\oplus d}$. \end{proof}

\begin{lemma}\label{l.-1}
In Definition \ref{d.sP}, let $\sU$ be the tautological vector bundle on ${\rm Gr}(d, T \sR)$.
Then the  restriction of the tautological vector bundle $\sU$  to a fiber of the $\BP^1$-bundle $\psi: \sP \to \sR$ is isomorphic to $\sO(-1) ^{\oplus d}$. \end{lemma}

\begin{proof}
Since $\sP$ in Definition \ref{d.sP}  arises from the setting of Lemma \ref{l.tensor} in  Corollary \ref{c.tensor}, the tautological vector bundle splits into a direct sum of tautological line bundles on a fiber of $\psi$ by Lemma \ref{l.tautological}.  \end{proof}

\begin{notation}\label{n.PK} In Definition \ref{d.sP},  fix a 1-dimensional vector space $\I$ with a fixed nonzero vector $\bi \in \I$ and define $ \V = \W \oplus \U \oplus \I$. Identify $\G \stackrel{\log}{=} \W \oplus \U$ with the affine open subset
$\BP \V \setminus \BP (\W \oplus \U)$ by sending $(x^W, x^U)$ to $[x^W: x^U: \bi] \in \BP \V$. \begin{itemize}
\item[(i)]
For each horizontal line $\ell^w_x \subset \G$, we have a uniquely determined projective line $\bl^w_x \subset \BP \V$ such that $\ell^w_x = \G \cap \bl^w_x$. Thus there exists an injective holomorphic map $\vartheta$ from the set $\sR$ of $S$-lines on $\G$ to $ {\rm Gr}(2, \V).$
Let $\sR' \subset {\rm Gr}(2, \V)$ be its image.
\item[(ii)]
Let $\sP' \subset \BP \V \times \sR'$ be the universal family of projective lines parameterized by $\sR'$ with the universal family morphisms $$ \sR' \stackrel{\rho'}{\longleftarrow} \sP' \stackrel{\mu'}{\longrightarrow}  \BP \V.$$
We have a natural injection $\sP \setminus \Sigma \subset \sP'$, which can be extended to a bijective holomorphic map   $\theta: \sP \to \sP'.$  In fact, when $\ell \in \sR$ is the point representing an $S$-line $\ell^w_x,$   the $\theta$-image of the point $\sigma (\ell) \in \sP_{\ell}$  is defined to be the point  $\bl^w_x  \cap \BP (\W \oplus \U) \subset \BP \V$.  This gives the following commutative diagram.
$$ \begin{array}{ccccc} \sR' & \stackrel{\rho'}{\longleftarrow} & \sP' & \stackrel{\mu'}{\longrightarrow}   & \BP \V \\
\vartheta \uparrow & & \uparrow \theta & &   \\  \sR & \stackrel{\psi}{\longleftarrow} & \sP &  & \cup \\
\parallel &  & \cup  &  &  \\
  \sR & \stackrel{\rho}{\longleftarrow} & \sS & \stackrel{\mu}{\longrightarrow} & \G. \end{array} $$
  \end{itemize}
\end{notation}

The following lemma says that the vector bundle $\sV^o$ on $\sP \setminus \Sigma$ is a meromorphic object on $\sP,$ namely, it does not have essential singularity along $\Sigma$.

\begin{lemma}\label{l.sV}
For a vector bundle $V$ on a complex manifold, denote by $\sO(V)$ the sheaf of local holomorphic sections of $V$.  In Definition \ref{d.sP}, there exists a saturated coherent subsheaf $\sV^{\rm sh} \subset \sO(T \sP)$ such that $\sV^{\rm sh}|_{\sP \setminus \Sigma} = \sO(\sV^o)$ and sections of $\sV^{\rm sh}$ are tangent to $\Sigma$. \end{lemma}

  \begin{proof}  The leaves of the foliation $\sV^o$ on $\sS$ agrees with the fibers of $\mu'\circ \theta$
  restricted to $\sS \subset \sP$ in Notation \ref{n.PK}. Thus we may define $\sV^{\rm sh} \subset \sO(T \sP)$ as the saturation of the sheaf of vector fields tangent to fibers of $\mu' \circ \theta$. \end{proof}

\begin{definition}\label{d.sJ}
In Definition \ref{d.sP}, define a vector subbundle $\sJ \subset T \sP$ of rank $d+1$ as follows.
For a point $[V] \in \sP$ corresponding to a $d$-dimensional subspace $V \subset T_{\ell} \sR$ with $\ell = \psi([V])$, the fiber $\sJ_{[V]} \subset T_{[V]} \sP$ is defined to be
$$\sJ_{[V]} := ({\rm d}_{[V]} \psi)^{-1}(V),$$ where ${\rm d}_{[V]} \psi: T_{[V]} \sP \to T_{\ell} \sR$ is the differential of $\psi$ at $[V]$.  Denoting by $T^{\psi} \subset T\sP$ the kernel of ${\rm d} \psi$, we have the exact sequence of vector bundles on $\sP$,
\begin{equation}\label{e.split}  0 \to T^{\psi} \to \sJ \to \sU|_{\sP} \to 0, \end{equation} where $\sU|_{\sP}$ is the restriction of the tautological vector bundle $\sU$
on ${\rm Gr}(d, T \sR)$ whose fiber $\sU_{[V]}$ at $[V] \in \sP$ is just $V \subset T_{\ell}\sR.$
\end{definition}

\begin{proposition}\label{p.codim2}
In Definition \ref{d.sP}, there exists a
 vector subbundle  $\sV \subset T \sP$ of rank $d$ such that \begin{itemize} \item[(i)]
 $\sV$ is contained in the vector bundle $\sJ$ in Definition \ref{d.sJ} and  splits the exact sequence (\ref{e.split}); \item[(ii)]
 the restriction $\sV|_{\sP \setminus \Sigma}$ agrees with $\sV^o$: and    \item[(iii)]  the restriction $\sV|_{\Sigma}$
 agrees with $\sT \subset  T \Sigma$. \end{itemize} \end{proposition}

\begin{proof}
As $\sV^{\rm sh}$ in Lemma \ref{l.sV} is saturated, there exists an analytic subset $Z \subset \sP$ of codimension $\geq 2$ and a vector subbundle $\sV' \subset T(\sP \setminus Z)$ such that \begin{equation}\label{e.sh} \sV^{\rm sh}|_{\sP \setminus Z} = \sO(\sV').\end{equation} By $\sV'|_{\sP \setminus \Sigma} = \sV^o,$
the distribution $\sV'$ is integrable. Since sections of $\sV^{\rm sh}$ are tangent to $\Sigma$, the leaves of $\sV'$ foliate $\Sigma \setminus Z \cong \sR \setminus \psi(Z)$.
For each point $\ell \in \sR \setminus \psi(Z)$, the family of submanifolds in $\sR$  through $\ell$ $$\{ \rho( \mu^{-1}(x_t)) \subset \sR \mid x_t \in \mu(\rho^{-1}(\ell)) \subset \G\}$$ are $\psi$-images of leaves of $\sV'$ on $\sP$. Thus the limit of this family gives the leaf of $\sV'$ through $\ell$ on $\sR \setminus \psi(Z)$.
By Corollary \ref{c.sT}, this limit is tangent to $\sT$. It follows  that  \begin{equation}\label{e.sV'} \sV'|_{\Sigma \setminus Z} = \sT|_{\Sigma \setminus Z}. \end{equation} In particular, the restriction $\sV'|_{\Sigma \setminus Z}$ splits (\ref{e.split}) on $\Sigma \setminus Z$. Note that  $\sV^o \subset \sJ|_{\sP \setminus \Sigma}$ splits (\ref{e.split}) on $\sP \setminus \Sigma$.  Thus $\sV'$ is a subbundle of $\sJ|_{\sP \setminus Z}$ and it splits (\ref{e.split}) on $\sP \setminus Z$.
Since $Z $ has codimension at least $2$ in $\sP$, the vector bundle $\sV'$    splitting (\ref{e.split}) on $\sP \setminus Z$ can be  extended  to  a vector subbundle $\sV \subset T \sP$
splitting (\ref{e.split}) on $\sP$  by Hartogs extension (see p. 409 of \cite{HM98} or Proposition 5 in \cite{Li}).   This proves (i). (ii) follows from $\sV|_{\sP \setminus Z} = \sV'$ and (\ref{e.sh}).   (iii) follows from (\ref{e.sV'}).   \end{proof}

\section{A partial compactification of $\G$ and the proof of  Theorems \ref{t.analytic} and \ref{t.algebraic}}\label{s.algebraic}

The following construction is similar to Definition \ref{d.sL}, replacing the group $\A_s$ by $\T_s$ introduced in Definition \ref{d.sT}. 

\begin{definition}\label{d.sH} In Definition \ref{d.sT},
let $\T \subset S \times \G$ be the family of subgroups $$\{\T_s \subset \{s\} \times \G \mid s \in S \}$$ and $\nu: \J \to S$ be the family of coset spaces. In other words, we have a sequence of morphisms $$\begin{array}{ccccc} \T & \subset & S \times \G & \longrightarrow & \J \\
\downarrow & & \downarrow & & \downarrow \nu \\ S & = & S & = & S\end{array}$$
which gives an exact sequence
$$0 \to  \T_s \longrightarrow \G \longrightarrow \G/\T_s = \nu^{-1}(s) \to 0$$ over each point $s \in S$.
 Then we can regard $\J$ as the  space of the leaves of the foliation $\sT$ on $\Sigma$ and obtain a natural surjective holomorphic map $\xi: \Sigma \to \J$ such that $\sT = {\rm Ker}({\rm d} \xi) \subset T \Sigma.$  \end{definition}

\begin{theorem}\label{t.quotient}
Using the terminology of  Proposition \ref{p.codim2} and Definition \ref{d.sH},
let $\X$ be the disjoint union $\G \cup \J$ and let $\widehat{\mu}:  \sP \to \X$ be the set-theoretical map
defined by the combination of $ \mu: \sP \setminus \Sigma \stackrel{\chi}{=} \sS \stackrel{\mu}{\to} \G$ and $\xi: \Sigma  \to \J:$
$$\begin{array}{ccccc}
\sP & = & (\sP \setminus \Sigma) & \cup & \Sigma \\
\widehat{\mu} \downarrow & & \mu \downarrow & & \downarrow \xi \\
\X & = & \G & \cup & \J.\end{array} $$
Then $\X$ has a natural  complex manifold structure such that  $\widehat{\mu}$ is holomorphic and $\widehat{\mu}$ sends every fiber of $\psi: \sP  \to \sR$  to a smooth unbendable rational curve in $\X$ with normal bundle isomorphic to $\sO(1)^{\oplus d} \oplus \sO^{\oplus (n-1-d)}$. Furthermore, the left translation of $\G$ can be extended to a left $\G$-action on $\X$.  \end{theorem}

\begin{proof}
By definition, the set $\X$ is the union of two complex manifolds $\G$ and $\J$. We need to patch the two complex manifolds to define a complex structure on $\X$.

Fix a point $\ell \in \sR$ and consider the fiber $\sP_{\ell} := \psi^{-1}(\ell)$.
As $\sV$ from Proposition \ref{p.codim2} splits the exact sequence (\ref{e.split}), we can find an open neighborhood  $O_{\ell} \subset \sP$ of $\sP_{\ell}$ such that leaves of the foliation given by $\sV$ define a submersion $\zeta_{\ell}: O_{\ell} \to U_{\ell}$ onto a complex manifold $U_{\ell}$. There is a natural set-theoretical injection $\epsilon_{\ell}: U_{\ell} \to \X$ which sends $\zeta_{\ell}(\sigma(\ell)) $ to $  \xi (\sigma(\ell)) \in \J  \subset \X$.   By transplanting the complex manifold structure of $U_{\ell} $ to the subset $\epsilon_{\ell}(U_{\ell})$, we obtain a neighborhood of $\xi (\sigma(\ell))$ in $\X$ with a complex manifold structure which is compatible with the complex structures of $\G$ and $\J$. The germ of this complex structure at $\xi (\sigma(\ell)) \in \X$ is independent of the choice of $O_{\ell}$.

We need to check that for two different points $\ell \neq \ell' \in \sR$ with $\xi (\sigma(\ell))= \xi (\sigma(\ell'))  \in \J$, the complex structure in a neighborhood of  $\xi (\sigma(\ell))= \xi (\sigma(\ell'))$ in $\X$  arising from $\ell$   coincides with the one arising from $\ell'$.
When $\ell \in \sR$ represents a horizontal  line of the form $\ell^w_x$ with $[w] \in S$, we can find an element $g \in \T_s$ such that $\ell'$ represent the line $\ell^w_{x \cdot g}.$ Note that there exists a natural $\G$-action on $\sP$ which preserves the distribution $\sV \subset T\sP$.
From $$\ell^w_{x \cdot g} = x \cdot g \cdot \ell^w_o =
(x \cdot g \cdot x^{-1}) \cdot \ell^w_x,$$ we see that $(x \cdot g \cdot x^{-1}) \cdot \sP_{\ell} = \sP_{\ell'}$ by the natural  $\G$-action on $\sP$.  Thus we may choose $O_{\ell'}$ as $(x \cdot g \cdot x^{-1}) \cdot O_{\ell}$ by the left action of $x \cdot g \cdot x^{-1} \in \G$. As the  $\G$-action  preserves the foliation $\sV$ on $\sP$, it is compatible with the map $\xi$ and the left action of $x \cdot g \cdot x^{-1}$ gives a biholomorphic map between the two complex structures arising from $\ell$ and $\ell'$. Thus we obtain a well-defined complex  structure on $\X$ compatible with those of $\G$ and $\J$.
The map $\widehat{\mu}$ is holomorphic by the definition of the complex structure on $\X$.

The image of $\sP_{\ell}$ in $U_{\ell}$ under $\zeta_{\ell}$ is a smooth rational curve $C_{\ell}.$
The union of deformations of $C_{\ell}$ in $\X$ through a fixed point $x \in \G \subset \X$
is an analytic subset of dimension at least $d+1$ because it includes the union of $S$-lines in $\G$ through $x$. It follows that the normal bundle $N_{C_{\ell}}$ of $C_{\ell} \subset \X$ contains at least $d$ positive factors. The normal bundle $N_{\sP_{\ell}}$ of $\sP_{\ell}$ in $\sP$ is a trivial vector bundle.  From the exact sequence
$$ 0 \to \sV|_{\sP_{\ell}} \to N_{\sP_{\ell}} \to \widehat{\mu}^* N_{C_{\ell}} \to 0$$ and the isomorphisms $$\sV|_{\sP_{\ell}} \cong \sU|_{\sP_{\ell}} \cong \sO(-1)^{\oplus d} $$ from Lemma \ref{l.-1} and Proposition \ref{p.codim2} (i), we conclude that
 $N_{C_{\ell}}$ is isomorphic to $\sO(1)^{\oplus d} \oplus \sO^{\oplus{(n-1-d)}}.$
Thus $\widehat{\mu}$ sends a fiber of $\psi$ to a smooth unbendable rational curve with normal bundle
 isomorphic to $\sO(1)^{\oplus d} \oplus \sO^{\oplus{(n-1-d)}}$.

 Finally, it is obvious that the left translation on $\G$ can be extended to a $\G$-action on $\X$.
\end{proof}

\begin{proof}[Proof of Theorem \ref{t.analytic}]
The proof is easy from Theorem
 \ref{t.quotient}. We are assuming that the submanifold $S \subset \BP \W$ is a projective variety and $\mu: \sS \to \G$ is a proper holomorphic map.  Set $\sX := \X$, $\sX^o := \G$ and $D := \sW$.
 Viewing $$ \sR \stackrel{\psi}{\longleftarrow} \sP \stackrel{\widehat{\mu}}{\longrightarrow} \X$$ as the universal family morphisms of a set of unbendable rational curves in $\X$, we have a natural inclusion $\sR \subset {\rm RatCurves}(\X)$. From $$ \dim \sR = n-1 + d = \dim H^0(\BP^1, \sO(1)^{\oplus d} \oplus \sO^{\oplus (n-1-d)}),$$ we see that $\sR$ is an open subset in an irreducible component $\sK$ of ${\rm RatCurves}(\X)$. Since $\widehat{\mu} |_{\chi (\sS)}$ coincides with the proper holomorphic map $\mu: \sS \to \G$, we see that $\sK$ is a family of minimal rational curves on $\X$ with the open set $O_{\sK} \subset \X$ of Definition \ref{d.mrc} (v) equal to $\G \subset \X$.
 It is clear that $\sK$ is subordinate to the distribution $D = \sW$ on $\G$.
 \end{proof}

\begin{proof}[Proof of Theorem \ref{t.algebraic} ]
We are assuming that the submanifold $S \subset \BP \W$ is a projective variety.  Then the set $\sR' \subset {\rm Gr}(2, \V)$ in Notation \ref{n.PK} is quasi-projective. Let $\bR$ be its closure in ${\rm Gr}(2, \V)$ and $\bP \subset \BP \V \times \bR$ be the universal family of projective lines parameterized by $\bR$ with the universal family morphisms $$ \bR \stackrel{\widetilde{\rho}}{\longleftarrow} \bP \stackrel{\widetilde{\mu}}{\longrightarrow}  \BP \V.$$ Then $\bP$ is a projective variety.
  The universal family of projective lines in $\BP \V $ parameterized by ${\rm Gr}(2, \V)$ can be naturally identified with $\BP(T \BP \V)$. Thus we have a natural inclusion $\bP \subset \BP (T \BP \V)$ satisfying the diagram $$ \begin{array}{ccccc}
   \sS & \subset & \bP & \subset & \BP T (\BP \V) \\
  \mu \downarrow & &  \downarrow \widetilde{\mu} & & \downarrow \\
  \G & \subset  & \BP \V & = & \BP \V. \end{array} $$ The fibers of $\widetilde{\mu}$ over $\G \subset \BP \V$ are isomorphic to the projective variety
  $S \subset \BP \W$. Let $\bY^o$  be the quasi-projective subvariety of  the Chow variety  of $\bP$ which parameterizes the fibers of $\mu: \sS \to \G$. There is a natural biregular isomorphism $\varsigma: \G \to \bY^o$ which sends $x \in \G$ to the point in $\bY^0$ corresponding to $\mu^{-1}(x) \subset \bP$.   Let $\bY$ be the normalization of the closure of $\bY^o$ in the Chow variety of $\bP$.

Let $\sV$ be the distribution on $\sP$ from Proposition \ref{p.codim2}.
Note that the Zariski-open subset $\bY^o$ of $\bY$ parameterizes the $d$-dimensional irreducible cycles in $\bP$ corresponding to the leaves of $\sV$  on $\sP \setminus \Sigma$ and
the leaves of $\sV$ on $\Sigma$ are leaves of $\sT$. It follows that  there are
a finite number of irreducible hypersurfaces $\bH_1, \bH_2, \ldots, \bH_m \subset \bY \setminus \bY^o$ with nonempty Zariski-open subsets $\bH_i^o \subset \bH_i, 1 \leq i \leq m$,  such that each point of  $\bH^o_i$ corresponds to  a   $d$-dimensional cycle  in $\bP$ of the  form
\begin{equation}\label{e.cycle} a_1 F_1 + \cdots + a_r F_r \end{equation} for some positive integers $r, a_1, \ldots, a_r$ and distinct $d$-dimensional  irreducible subvarieties $F_1, \ldots, F_r$ in $\bP$  such that $F_1$ contains a leaf $F^o_1$ of
 $\sT$ on $\Sigma \subset \bP$ as a Zariski-open subset.
  Furthermore, we can assume that a  general leaf $F_1^o$  of $\sT$  appears as a Zariski-open subset of $F_1$ in the expression (\ref{e.cycle}) of a member of $\bH^o_i$ for each $1 \leq i \leq m$.

Consider the complex manifold $\X$ from Theorem  \ref{t.quotient}.
Regard the natural inclusion $\varepsilon:  \G \subset \X$ as a holomorphic map defined on the open subset $\G \stackrel{\varsigma}{=} \bY^o  \subset \bY$.  We can extend the holomorphic map $\varepsilon: \bY^o \to \X$ to a holomorphic map $$\widetilde{\varepsilon}: \bY^o \cup \bH_1^o \cup \cdots \cup \bH_m^o \to \X$$  by sending a point of $\bH^o_i$  corresponding to a cycle of the form (\ref{e.cycle}) to the point of $\J  \subset \X$ corresponding to the leaf $F^o_1$ of $\sT$.
This map $\widetilde{\varepsilon}$  is a continuous extension of $\varepsilon$, hence it is a holomorphic extension of $\varepsilon$ because $\bY$ is normal.
By our assumption that   a general leaf of $\sT$  appears as $F^o_1$ for a member  of $\bH^o_i$ for each $1 \leq i \leq m,$ the image of $\widetilde{\varepsilon}$ is $\X \setminus Z$ for some analytic set $Z \subset \J$ of codimension at least 2 in $\X$. After replacing $\bH^o_i$ by a smaller Zariski-open subset in $\bH_i$ if necessary, we can assume that $$\widetilde{\varepsilon}:
 \bY^o \cup \bH_1^o \cup \cdots \cup \bH_m^o \to \X \setminus Z$$ is a biholomorphic map by Zariski Main Theorem (e.g. Corollary in Section 4.9 of \cite{Fi}).

Let $\bX \to \bY$ be a desingularization of $\bY$ which is isomorphic over the smooth locus of $\bY$.
Let $\bX^o \subset \bX$ be the Zariski-open subset corresponding to $\widetilde{\varepsilon}^{-1}(\G)$ and $D \subset T\bX^o$ be the distribution corresponding to $\sW \subset T\G$.
The unbendable rational curves in $\X$ given by general fibers of $\psi: \sP \to \sR$ in Theorem \ref{t.quotient}  are disjoint from $Z$ and their images under $\widetilde{\varepsilon}^{-1}$ give unbendable rational curves on $\bX$. They determine a  family $\sK$ of unbendable rational curves on $\bX$ which is subordinate to the distribution $D$ and the VMRT at a point of $\bX^o$ is isomorphic to $S \subset \BP \W$.
\end{proof}

\bigskip
{\bf Acknowledgment}
I would like to thank Qifeng Li for valuable comments on the first draft of the paper, especially pointing out the equivalence of (\ref{e.HM}) and (\ref{e.Hw}).

\bigskip
Institute for Basic Science

Center for Complex Geometry

Daejeon, 34126

Republic of Korea

jmhwang@ibs.re.kr
\end{document}